\documentclass{amsart}
\usepackage{amssymb}
\usepackage{url}

\newcommand{\be}{\begin}
\newcommand{\en}{\end}
\newcommand{\bpr}{\begin{proof}}
\newcommand{\epr}{\end{proof}}

 \newtheorem{theorem}{Theorem}

 \newtheorem{example}[theorem]{Example}

\newtheorem{pro}[theorem]{Proposition}

 \newcommand{\Z}{{\mathbb Z}}
 \newcommand{\C}{{\mathbb C}}
 \newcommand{\N}{{\mathbb N}}
 \newcommand{\T}{{\mathcal T}}
\newcommand{\K}{{\mathbb K}}
\newcommand{\FT}{{\cal F\cal T}}
\newcommand{\cal}{\mathcal}
\newcommand{\X}{\mathcal X}
\renewcommand{\S}{\mathcal S}

\newcommand{\ve}{\varepsilon}

\newcommand{\B}{\cal B}

\newcommand{\D}{\cal D}
\newcommand{\pf}{P\hspace*{-.02in}f}

\newcommand{\la}{\left\langle}
\newcommand{\ra}{\right\rangle}

\newcommand{\comment}[1]{}

\title{Generating sets for coordinate rings of character varieties}
\author{Adam S. Sikora}


\date{}

\begin{document}
\thispagestyle{empty}

\begin{abstract}
We find finite, reasonably small, generator sets of the coordinate rings of $G$-character varieties of finitely generated groups for all classical matrix groups $G$. This result together with the method
of Gr\"obner basis gives an algorithm for describing character varieties by explicit polynomial equations.

Additionally, we describe finite sets of generators of the fields of rational functions on $G$-character varieties for all exceptional algebraic groups $G.$
\end{abstract}

\maketitle

%
\section{Character Varieties}
\label{s_char_var}
%

Throughout the paper $G$ will be an affine reductive group over an algebraically closed field $\K$ of characteristic zero and $\Gamma$ will be a (discrete) group generated by $\gamma_1,...,\gamma_N.$
The space of all $G$-representations of $\Gamma$ forms an algebraic subset, \\ $Hom(\Gamma,G),$ of $G^N$, called the $G$-representation variety of $\Gamma$.
The group $G$ acts on this set by conjugating representations and the categorical quotient of that action $$X_G(\Gamma)=Hom(\Gamma,G)//G$$ is the $G$-character variety of $\Gamma;$ cf. \cite{S-char} and the references within.

Due to the ubiquity of applications of character varieties in low-dimensional topology, geometry, gauge theory, and quantum field theories
one is interested in an explicit description of them by polynomial equations, or, equivalently a description of $\K[X_G(\Gamma)]$ by generators and relations.

In this paper, we describe generating sets of the coordinate rings $\K[X_G(\Gamma)]$ for all classical matrix groups $G$.
(A discussion of generating sets for finite quotients of classical matrix groups appears in \cite{S-quot}.)
Additionally, we describe generators of the fields of rational functions on $G$-character varieties for all exceptional groups $G.$

We do not discuss here finding relations between generators. An algorithmic solution to this problem is given by the theory of Gr\"obner bases. (However, due to its computational complexity, it is difficult to apply this method in practice for large generating sets.) Character varieties of abelian groups are studied in more detail in \cite{S-abel}.

Let $T$ be a maximal torus in $G$.

\begin{example}\textup{(\cite[6.4]{St}, \cite[Example 42]{S-char})}\label{Z}
If $G$ is connected, then ${\cal X}_G(\Z)=T/W,$ where $T/W$ is the quotient of $T$ by the action of the Weyl group $W$ of $G$ with respect to $T$. Consequently, $\K[X_G(\Z)]$ is isomorphic to the representation ring of $G$.
\end{example}

For more concrete algebraic description of $X_G(\Z)$, we refer the reader to \cite[Thm 3.12 and Eg. 3.15]{Po}.

For a matrix group $G\subset GL(n,\K)$ and $\gamma\in \Gamma$, we define $\tau_\gamma: X_G(\Gamma)\to \K$
by $\tau_\gamma(\rho)=tr(\rho(\gamma)).$

The following is well known; see for example \cite{Sib}, \cite{BH},\cite[Thm 8.4]{PS}, and a new proof in Sec. \ref{s_proof_SL2}. (A larger generating set was found in \cite{Ho}.)

\be{pro}\label{SL2-free} Let $\Gamma$ be generated by $\gamma_1,...,\gamma_N.$\\
\textup{(1)} The $\K$-algebra $\K[X_{SL(2,\K)}(\Gamma)]$ is generated by $\tau_{\gamma_i},$ for $i=1,..,N,$ $\tau_{\gamma_i\gamma_j}$ for $i<j$ and $\tau_{\gamma_i\gamma_j\gamma_k},$
for $i<j<k.$\\
\textup{(2)} If $\Gamma$ is abelian, then $\K[X_{SL(2,\K)}(\Gamma)]$ is generated by $\tau_{\gamma_i},$ for $i=1,..,N,$ and $\tau_{\gamma_i\gamma_j}$ for $i<j.$
\en{pro}

Furthermore, the above generating sets are minimal, cf. \cite[Thm 8.6]{PS}.
Minimal generating sets for $G$-character varieties of free groups are studied in \cite{ADS,AP,BD, La,Na,T1} for $G=SL(3,\K)$ and $GL(3,\K)$ and in \cite{Dj, DL, DS, T2} for $G=GL(4,\K).$

\noindent {\bf Acknowledgments}
We would like to thank S. Lawton for helpful conversations.
He has discovered some of the results of the next section independently of us.

%
\section{Generators of $\K[X_G(\Gamma)]$ for Classical Groups -- Main Results}
\label{s_main}
%

The nil index of an associative ring $R$ (without identity) is the smallest positive integer $m$ such that $r^m=0$ for all $r\in R$. The nilpotency index of $R$, $\nu(R),$ is the smallest positive integer $n$ such that the product of every $n$ elements of $R$ vanishes.

Let $\nu_n$ be the lowest upper bound on the nilpotency index of associative $\K$-algebras of nil index $n.$
By \cite{Ra} (page 759 in the English translation), $$\nu_n\leq n^2.$$ (The Nagata-Higman theorem, established earlier, implies the finiteness of $\nu_n$.) Furthermore, Kuzmin's conjecture, stating that
\be{equation}\label{e-Ku}
\nu_n= \frac{(n+1)n}{2}
\en{equation}
for all $n$, \cite{Ku}, holds for $n\leq 4$, \cite{Du,VL}.

\begin{theorem}\label{basic}
\textup{(1)} Let $\B$ denote the set all elements of $\Gamma$ which are represented by words in $\gamma_1,...,\gamma_N$ (without negative exponents) of length at most $\nu_n$. Let $\B'$ be a subset of $\B$ which contains representatives of all conjugacy classes (in $G$) of elements of $\B$.
Then $\K[X_{SL(n,\K)}(\Gamma)]$ is generated by $\tau_{\gamma},$ for $\gamma\in \B'.$\\
\textup{(2)} For abelian $\Gamma$, $\K[X_{SL(n,\K)}(\Gamma)]$ is generated by $\tau_\gamma$ for $\gamma$ represented by words (without negative exponents) of length at most $n$.
\end{theorem}

This theorem is a direct consequence of Theorem \ref{sln-main}, which provides more efficient generating sets.

\be{rem}\label{rem_gln}
If $\Gamma$ is generated by $\gamma_1,...,\gamma_N$ and
$\K[X_{SL(n,\K)}(\Gamma)]$ is generated by
$\tau_{\gamma},$ for $\gamma$ in some set $\B\subset \Gamma,$
then $\K[X_{GL(n,\K)}(\Gamma)]$ is generated by
$\tau_{\gamma},$ for $\gamma\in \B,$ and by the functions $\sigma_i(\rho)=det(\rho(\gamma_i))^{-1},$ $i=1,...,N.$
\en{rem}

Recall that
$$O(n,\K)=\{A\in GL(n,\K): AA^T=I\},\quad SO(n,\K)=O(n,\K)\cap SL(n,\K),$$
and
$$Sp(n,\K)=\{A\in GL(n,\K): AJA^T=J\}$$
where $n$ is even and $J$ is a non-degenerate skew-symmetric matrix, eg.\\
$J=\left( \begin{array}{cc} 0 & I_{n/2}\\ -I_{n/2} & 0\\ \end{array}\right)$.

We will call $SO(n,\K)$ an {\em even} or {\em odd} special orthogonal group, depending on the parity of $n$.

For any matrix group $G\subset GL(n,\K),$ let $\T_G(\Gamma)\subset \K[X_G(\Gamma)]$ be the
subalgebra generated by the functions $\tau_\gamma,$ for all $\gamma\in \Gamma.$
(Keep in mind that $\T_G(\Gamma)$ depends on the specific embedding $G\hookrightarrow GL(n,\K).$) We call  $\T_G(\Gamma)$ the {\em $G$-trace algebra of $\Gamma$.}

By Theorem \ref{basic}, $\K[X_{SL(n,\K)}(\Gamma)=\T_{SL(n,\K)}(\Gamma)$ for every $\Gamma.$

\be{theorem}\textup{(\cite[Thm A.1]{FL}, see also Sec. \ref{spso-proof})}\label{spso}
If $G$ is symplectic, orthogonal, or odd special orthogonal, then $\K[X_G(\Gamma)]=\T_G(\Gamma)$ for every $\Gamma.$
\en{theorem}

\be{rem} 
It is easy to see that
$\T_{GL(n,\K)}(\Z)=\K[x_1,...,x_n]^{S_n}$, while\\
$\K[X_{GL(n,\K)}(\Gamma)]=\K[x_1^{\pm 1},...,x_n^{\pm 1}]^{S_n}$ by Example \ref{Z}. Hence $\T_{GL(n,\K)}(\Z)$ is a proper subalgebra of $\K[X_{GL(n,\K)}(\Z)].$
\en{rem}

It is easy to see that the trace algebras satisfy the following nice property: every embedding $G\subset G'\subset GL(n,\K),$ induces an epimorphism $\T_{G'}(\Gamma)\to \T_{G}(\Gamma)$.
(The induced map $\K[X_{G'}(\Gamma)]\to \K[X_{G}(\Gamma)]$ does not have to be an epimorphism. That happens for example for $G=SO(2,\K)\subset SL(2,\K)=G'$ -- see the example following Theorem \ref{basic_SO(2m)}.)




\be{corollary}\label{trace-epi}
Any generating set of $\K[X_{SL(n,\K)}(\Gamma)]$ yields a generating set of $\K[X_{G}(\Gamma)]$ for $G$ symplectic and odd special orthogonal groups through the epimorphism
$$\K[X_{SL(n,\K)}(\Gamma)]= \T_{SL(n,\K)}(\Gamma)\to \T_G(\Gamma)=\K[X_G(\Gamma)].$$
\en{corollary}

Smaller generating sets for symplectic groups are given by Theorem \ref{sp-main} in Section \ref{s_spso}.

The generators of $\K[X_G(\Gamma)]$ for even special orthogonal groups are more difficult to describe, since unlike for other classical groups, $\T_G(\Gamma)$ is usually a proper subalgebra of $\K[X_G(\Gamma)]$ in that case.

For an even $n$, consider a function $Q_n: M(n,\K)^{n/2}\to \K,$ on the Cartesian product of $n/2$ copies of $n\times n$ matrix algebras, given for matrices $A_1,...,A_{n/2}$ by\vspace*{.1in}

$Q_n(A_1,...,A_{n/2})=\sum_{\sigma\in S_n} sn(\sigma)
 (A_{1,\sigma(1),\sigma(2)}-A_{1,\sigma(2),\sigma(1)})...\hspace*{.2in}$\vspace*{-.2in}\\
\be{equation}\label{def-Q}
\hspace*{2.8in} (A_{n/2,\sigma(n-1),\sigma(n)}-A_{n/2,\sigma(n),\sigma(n-1)}),
\en{equation}
where $A_{i,j,k}$ is the $(j,k)$-th entry of $A_i$ and $sn(\sigma)=\pm 1$ is the sign of $\sigma.$

Since $Hom(\Gamma,SO(n,\K))\to \K$ sending $\rho$ to $Q_n(\rho(\gamma_1),...,\rho(\gamma_{n/2}))$ is a regular function invariant under the conjugation of $\rho$ by elements of $SO(n,\K)$, it factors to a function
on $X_{SO(n,\K)}(\Gamma)$ which we denote by $Q_n(\gamma_1,...,\gamma_{n/2}).$
The following is a consequence of our Theorem \ref{so2n-main}.

\begin{theorem}\label{basic_SO(2m)}
For $n$ even, $\K[X_{SO(n,\K)}(\Gamma)]$ is generated as a $\T_{SO(n,\K)}(\Gamma)$-algebra by
$Q_n(w_1,...,w_{n/2}),$ for all $w_1,...,w_{n/2}\in \Gamma$ represented by words of length at most $\nu_n-1$ in which the number of inverses is not larger than half the length of the word.
\end{theorem}

Consequently, $\K[X_{SO(n,\K)}(\Gamma)]$ is generated by $\tau_\gamma$, for $\gamma$ in a set $\B'$ described in Theorem \ref{basic} and by the functions $Q_n(w_1,...,w_{n/2})$ as above.

Elements $\tau_\gamma$ do not generate $\K[X_{SO(n,\K)}(\Gamma)]$ alone.
As pointed out to us by S. Lawton, that phenomenon appears even for $n=2$ and $\Gamma=\Z:$
Let $i\in \K$ be a primitive $4$th root of $1$ and let $\rho,\rho':\Z\to SO(2,\K)$ send $1$ to
$$A=\frac{1}{2}\left(\be{array}{cc} x+x^{-1} & i(x-x^{-1})\\ -i(x-x^{-1}) & x+x^{-1}\en{array}\right)\in SO(2,\K)$$
and to $A^T$, respectively, for some $x\ne \pm 1.$
Since $Q_2(\rho)=4(x-x^{-1})\ne Q_2(\rho')$, $[\rho]$ and $[\rho']$
are distinct points of $X_{SO(2)}(\Z).$ However $[\rho]$ and $[\rho']$
are not distinguished by $\tau_{\gamma}$ for any $\gamma\in \Z.$

For every representation $\phi:G\to GL(n,\K)$ and for every $\gamma\in \Gamma$, consider the regular function $\tau_{\gamma,\phi}: X_G(\Gamma)\to \K,$
$$\tau_{\gamma,\phi}([\rho])=tr(\phi\rho(\gamma))$$
and denote the $\K$-subalgebra of $\K[X_G(\Gamma)]$ generated by $\tau_{\gamma,\phi}$ over all $\gamma\in \Gamma$ and all representations $\phi$ of $G$ by $\FT_G(\Gamma).$
We call it the {\em full $G$-trace algebra of $\Gamma.$}
For a matrix group $\phi:G\hookrightarrow GL(n,\K)$, we have $\tau_{\gamma,\phi}=\tau_\gamma$ and
${\cal T}_{G}(\Gamma)\subset\FT_{G}(\Gamma).$
We prove in \cite{S-quot} that $\FT_{SO(4,\K)}(\Gamma)$ is a proper subalgebra of $\K[X_{SO(4,\K)}(\Gamma)]$ for every free group $\Gamma$ of rank $\geq 2.$
Consequently, the generators $Q_n(w_1,...,w_{n/2})$ of Theorem \ref{basic_SO(2m)} are not redundant in that case.
The following remains open:

\be{problem}\label{not_generated} Find explicit presentations (by generators and relations) of the extensions ${\cal T}_{SO(n,\K)}(\Gamma)\subset\FT_{SO(n,\K)}(\Gamma)\subset \K[X_{SO(n,\K)}(\Gamma)]$ for all $n$ and free groups $\Gamma$ of every rank.
\en{problem}

We continue the discussion of the relation between coordinate rings and the trace algebras in Sec. \ref{s_trace_a}.


Denote the Pfaffian of a matrix $X$ by $Pf(X).$ Here is an alternative description of $Q_n$ ($n$ even):

\be{pro}
\textup{(1)} $Q_n$ is a multi-linear, symmetric function such that
$$Q_n(X,...,X)=2^{n/2}(n/2)!\cdot Pf(X-X^T).$$
\textup{(2)} A function with these properties is unique, since it is the ``full polarization" of $2^{n/2}(n/2)!\pf(X-X^T).$\\
\textup{(3)} We have
$$\pf(X-X^T)=i^{n/2}(\tau_{D_+}(X)-\tau_{D_-}(X)),$$
where $\tau_{D_\pm}(X)$ denotes the trace of the image of $X\in SO(n,\K)$ under the $D_\pm$ representation of $SO(n,\K).$ It is a representation whose highest weight is twice that of $\pm$-half spin representation, cf. \cite[23.2]{FH}.
\en{pro}

\bpr (1) Since $Q_n(X,...,X)$ and $\pf(X-X^T)$ are conjugation invariant regular functions on $SO(n,\K)$, it is enough (by Example \ref{Z}) to prove their equality for $X$ in a maximal torus of $SO(n,\K)$, which is composed of block diagonal matrices, with diagonal blocks $$A_j=\frac{1}{2}\left(\be{array}{cc} x_j+x_j^{-1} & i(x_j-x_j^{-1})\\ -i(x_j-x_j^{-1}) & x_j+x_j^{-1}\en{array}\right),$$ for $j=1,...,n/2.$
Now
\be{equation}\label{qpf}
Q_n(X,...,X)=(2i)^{n/2}(n/2)! \prod_{j=1}^{n/2} (x_j-x_j^{-1}) =2^{n/2}(n/2)! \pf(X-X^T),
\en{equation}
by direct computation.

(2) Note that the function\\
$$Q_n(Y,Y,...,Y)=2^{n/2}(n/2)! Pf(Y-Y^T),$$ where $Y=\sum_{j=1}^{n/2} \alpha_jX_j,$
is a polynomial in variables $\alpha_1,...,\alpha_{n/2}$ and that\\ $Q_n(X_1,...,X_{n/2})$ is $\frac{1}{(n/2)!}$ of the coefficient of $\alpha_1\cdot ...\cdot \alpha_{n/2}$.

(3) It is again enough to prove $\pf(X-X^T)=\tau_{D_+}(X)-\tau_{D_-}(X)$
for $X$ in the maximal torus. $\tau_{D_\pm}(X)$ is given by the character of $D_\pm$. By the formula above Corollary 7.8 in \cite{Ad},
$$\tau_{D_\ve}(X)=\sum x_1^{\ve_1}...x_{n/2}^{\ve_{n/2}}+\sigma,$$
where the sum is over all $\ve_1,...,\ve_{n/2}\in \{+1,-1\}$ such that
the sign of $\ve_1\cdot ...\cdot \ve_{n/2}$ coincides with $\ve$. Here, $\sigma$ is a sum of terms of smaller degree, which are the same for both $\ve$.
Therefore,
$$\tau_{D_+}(X)-\tau_{D_-}(X)=\prod_{j=1}^{n/2} (x_j-x_j^{-1})$$
and the statement follows now from equation (\ref{qpf}).
\epr

\section{Relation of $\K[X_G(\Gamma)]$ to the trace algebra}
\label{s_trace_a}

Finding finite generating sets of $\K[X_G(\Gamma)]$ for $G$ other than the groups discussed above is difficult
since the invariant theory of such groups is not fully understood. The following open problem is particularly important:

\be{que} Which algebraic reductive groups $G$ can be realized as matrix groups such that $\K[X_G(\Gamma)]=\T_G(\Gamma)$ for free groups $\Gamma$?
\en{que}

If the answer to this question is positive for $G\subset SL(n,\K)$ then the generators of Theorem \ref{basic} map onto the generators of $\K[X_G(\Gamma)]$
under the epimorphism $$\K[X_{SL(n,\K)}(\Gamma)]=\T_{SL(n,\K)}(\Gamma)\to \T_{G}(\Gamma)=\K[X_{G}(\Gamma)].$$

Our next result is relevant to the above question. We say that representations $\phi,\phi':\Gamma\to G$ are {\em twins} if there is an automorphism $\tau$ of $G$ such that $\phi'=\phi\tau.$ We say that these twins are {\em identical twins} if they are equivalent in $X_G(\Gamma)$.

\be{thm}[Proof in Sec. \ref{s_norm}]\label{normalization}
Let $G$ be a connected reductive group and let $\Gamma$ be a non-trivial free group.\\
\textup{(1)} If a representation $\rho: G\to SL(n,\K)$ is faithful and all twins of $\rho$ are identical, then
the natural embedding $\T_{\rho(G)}(\Gamma)\hookrightarrow \K[X_G(\Gamma)]$ is the integral closure of $\T_{\rho(G)}(\Gamma).$\\
\textup{(2)} If the kernel of $\rho: G\to SL(n,\K)$ contains a non-trivial central element, then $\K[X_G(\Gamma)]$ is a non-trivial extension of the integral closure of the trace algebra $\T_{\rho(G)}(\Gamma).$ (Note that if the Lie algebra of $G$ is simple, then the kernel of every non-faithful representation contains a non-trivial central element.)
\en{thm}

The condition of $\K[X_G(\Gamma)]$ being the integral closure of the trace algebra $\T_{\rho(G)}(\Gamma)$
is equivalent to $X_G(\Gamma)$ being the normalization of the algebraic variety\\
$Spec(\T_{\rho(G)}(\Gamma))$.

\be{exa}\label{exceptional}
 Let $G$ be one of the algebraic groups $E_7,E_8,F_4,G_2$ (i.e. any exceptional algebraic Lie group except for $E_6$).
 Then $Out(G)$ is trivial, cf. \cite[Sec 16.3]{Sp}, \cite[Sec. 2.1]{AdC}. Since $G$ is simple, every non-trivial representation $\rho$ of $G$ is faithful and $\K[X_G(\Gamma)]$ is the integral closure of the trace algebra $\T_{\rho(G)}(\Gamma)$ for such $\rho.$ Each of these groups has a unique representation $\rho$ of minimal dimension. It remains an open question whether $\T_{\rho(G)}(\Gamma)\hookrightarrow \K[X_G(\Gamma)]$ is an isomorphism for such $\rho$'s.
\en{exa}

\be{pro}[Proof in Sec. \ref{s_norm}]\label{E_6}
Let $G$ be the exceptional group $E_6$ and $\rho:G\to GL(n,\K)$ be a non-trivial representation of it of minimum dimension. (There are two non-equivalent representation like that!) Then $\K[X_G(\Gamma)]$ is the integral closure of $\T_{\rho(G)}(\Gamma)$.
\en{pro}

\section{Generators of $\K(X_G(\Gamma))$}
\label{s_rational}

If $G$ acts on an algebraic set $X$, then the field of fractions, $Q(\K[X]^G)$, of $\K[X]^G$ embeds into the field of $G$-invariant rational functions on $X,$ $\K(X)^G.$ Although that embedding is often proper, in case of character varieties we have:

\be{pro} For every connected reductive $G$ and a free group $\Gamma$, the embedding $\K(X_G(\Gamma))\hookrightarrow \K(Hom(\Gamma,G))^G$ is an isomorphism.
\en{pro}

\be{proof} If $\Gamma=\Z$ then $X_G(\Gamma)=T/W,$ where $T$ is a maximal torus in $G,$ cf. Example \ref{Z}.
The elements of $T/W$ distinguish orbits of the $G$-action on $Hom(\Gamma,G)$ in general position, cf. \cite[Sec II.2.1]{PV}. Now the statement follows from \cite[Prop II.3.4]{PV}.

For higher rank $\Gamma$ the proof is similar.
The set of irreducible representations in $Hom(\Gamma,G)$ is open, cf. \cite[Prop. 27]{S-char}. The elements of $\K[X_G(\Gamma)]$ distinguish irreducible representations, cf. \cite[Sec. 11]{S-char} and, hence,
they distinguish orbits of the $G$-action on $Hom(\Gamma,G)$ in general position as well.
\en{proof}

\be{thm}
Let $G$ be an exceptional group and let $\rho:G\to GL(n,\K)$ be a non-trivial representation. If $G=E_6$, then assume additionally that $\rho$ is of minimal dimension. If $\Gamma$ is free on $\gamma_1,...,\gamma_N,$ then
$\K(X_G(\Gamma))$ is generated by elements $\tau_{\rho,g}$ for $g$ represented by words without negative exponents of length $\leq \nu_n.$
\en{thm}

\be{proof}
By Example \ref{exceptional} and Proposition \ref{E_6}, $\K(X_G(\Gamma))$ is generated by the elements of $\T_{\rho(G)}(\Gamma)$.
We have $\rho(G)\subset SL(n,\K).$ (Indeed, since $G$ is simple and the kernel of $G \xrightarrow{\rho} GL(n,\K) \xrightarrow{det} \K^*$ cannot be $G$, it must be trivial.) Hence, $\T_{\rho(G)}(\Gamma)$ is an epimorphic image of $\K[X_{SL(n,\K)}(\Gamma)]$. Now the statement follows from Theorem \ref{basic}.
\en{proof}

The following is open, except for $G=GL(n,\K)$, cf. \cite{BK}.

\be{problem}
Find a transcendence basis of $\K(X_G(\Gamma))$ for all free groups $\Gamma$ and all reductive groups $G$.
\en{problem}

%
\section{Efficient generators for $SL(n)$-Character Varieties}
\label{s_sln}
%

Let $R(\Gamma, G)$ be the universal representation algebra of $\Gamma$ into $G$ and let $\rho_U:\Gamma\to G(R(\Gamma,G))$ be the universal representation, cf. \cite[Sec. 5]{S-char} (and \cite{S-sln,BH, LM} for $G=SL(n,\K)$). Here $G(R)$ denotes the $R$-points of the affine group scheme canonically associated with $G$.

There is a natural $G$-action on $R(\Gamma,G),$ c.f. \cite{LM, S-sln}.
Indeed, by the universality of $\rho_U,$ for every $g\in G$, $g\rho_U g^{-1}=G(f_g)\rho_U,$ for some homomorphism of $\K$-algebras $f_g: R(\Gamma,G)\to R(\Gamma,G).$ It is easy to see that $g\to f_g$ defines an algebraic $G$-action.

As shown in \cite{S-char}, there is the natural isomorphism $R(\Gamma,G)/\sqrt{0}\to \K[Hom(\Gamma,G)]$ restricting to the isomorphism $R(\Gamma,G)^G/\sqrt{0}\to \K[Hom(\Gamma,G)]^G=\K[X_G(\Gamma)]$.
Therefore,
$$\X_G(\Gamma)=Spec\, R(\Gamma,G)^G$$
is an affine algebraic scheme closely related with $X_G(\Gamma).$ The algebra of regular functions on it, $\K[\X_G(\Gamma)]$, is $R(\Gamma,G)^G.$

Assume that $G$ is a matrix group. For $\gamma\in \Gamma$, denote $Tr(\rho_U(\gamma))$ by $\sigma_\gamma.$ Since
$$f_g(Tr(\rho_U(\gamma)))=Tr\, G(f_g)\rho_U(\gamma)=Tr(g\rho_U(\gamma) g^{-1})=Tr(\rho_U(\gamma)),$$
for every $g\in G$, we have $\sigma_\gamma\in R(\Gamma,G)^G.$
Note that the natural projection $R(\Gamma,G)^G\to \K[Hom(\Gamma,G)]^G$ maps $\sigma_\gamma$ to $\tau_\gamma.$

Theorem \ref{sln-main} describes efficient generating sets of $\K[\X_G(\Gamma)]$ formed by elements $\sigma_\gamma$ for some $\gamma\in \Gamma.$ This theorem is a generalization of Theorem \ref{basic}.


We say that a semigroup $\S$ is graded if there is a homomorphism of semigroups
$$deg: \S\to \N=\{1,2,...\}.$$ We call it the {\em degree map}. (Note that no graded semigroup contains an identity.)

Let $\S$ be a graded semigroup and let $\S_d=\{\gamma\in \S: deg(\gamma)\leq d\}.$
Let $I$ be the two-sided ideal in $\K\S$ generated by elements
$z^n$ for all $z\in \K\S.$ Let $\D\subset \S$ be such that $\D\cap \S_d$ spans $\K\S_d/I\cap \K\S_d$ for all $d\in \N.$

\begin{theorem}[Proof in Sec. \ref{s_sln-main}]\label{sln-main}\ \\
Let $\S$ be a graded semigroup generated by elements $s_1,...,s_N$ of degree $1$ and
let $\eta: \S\to \Gamma$ be a homomorphism of semigroups whose image generates $\Gamma$
(as a group).\\
\textup{(1)} $\K[\X_{SL(n,\K)}(\Gamma)]$ is generated by $\sigma_{\gamma},$ for $\gamma\in \eta(\B)$, where $$\B=\{s_1,...,s_N\}\cup \D \cdot \{s_1,...,s_N\},$$
for any $\D$ as above. (Here, $A\cdot B$ means $\{a\cdot b: a\in A, b\in B\}.$)\\
\textup{(2)}(Stronger version) Let $\S_{d,r}$ be the set of all elements of $\S$ which are represented by words in $\gamma_1,...,\gamma_r$ of degree $\leq d.$ If $\D\cap \S_{d,r}$ spans
$\K\S_{d,r}/I \cap \K\S_{d,r}$ for all $1\leq r\leq N$ and all $d\geq 1$,
then it is enough to take $$\B=\{s_1,...,s_N\}\cup \bigcup_{r=1}^N \D_r \cdot s_r,$$
where $\D_r\subset \D$ is the set of the elements of $\D$ which are words in $s_1,...,s_r$.
\end{theorem}

Obviously, one can further reduce the set of generators by eliminating words in $\B$ which are related to others by cyclic permutations of letters.

For a group $\Gamma$ generated by $\gamma_1,...,\gamma_N$, one can take $\S$ to be the free semigroup on
$s_1,...,s_{N}$ and make $\eta$ send $s_1,...,s_N$ to the generators of $\Gamma$. However, for many $\Gamma$ one can take $\S$ to be a proper quotient of that free semigroup, which will result in a smaller set of generators of $\K[X_{SL(n,\K)}(\Gamma)].$

By taking $\D$ composed of all elements represented by words of length $< \nu_n$ and by considering the natural projection $\K[\X_{SL(n,\K)}(\Gamma)]\to \K[X_{SL(n,\K)}(\Gamma)]$ sending $\sigma_\gamma\to \tau_\gamma$ we get Theorem \ref{basic}(1). Part (2) of that theorem follows from the following:

\be{lem}\label{abel} If $S$ is commutative, then the nilpotency index of $\K\S/I$ is at most $n.$
\en{lem}

\bpr For any $x_1,...,x_n\in \S$, consider the polynomial $p:\K^n\to \K \S/I$, $$p(c_1,...,c_n)=(c_1x_1+...+c_nx_n)^n=0.$$
Since $\K$ is infinite, all monomials of $p$ vanish, including
$n!x_1...x_nc_1...c_n.$ Since $\K$ is of zero characteristic, $x_1...x_n=0.$
\epr

A specific $\D$ for a free semigroup $\S$ and $n=3$ is given in \cite[Prop 2]{Lo}.
For a free semigroup $\S$ and higher $n$, $\D$ and $\D_r$ can be found using the following lemma:

\be{lem}\label{algor}
For every $d\in\N,$ $I\cap \K\S_d$ is spanned by elements
$$v_0\left(\sum_{\sigma\in S_k} w_{\sigma(1)}^{\alpha_1}... w_{\sigma(k)}^{\alpha_k} \right)v_1,$$ for all $v_0,v_1,w_1,...,w_k\in S$, where $w_1,...,w_k$ are distinct, $\alpha_1\geq \alpha_2 \geq ...\geq \alpha_k\geq 1$, and $v_0+v_1+\sum_{i=1}^k \alpha_i deg(w_i)=d.$ (The words $v_0,v_1$ can be empty.)
\en{lem}

\bpr
Every element of $I$ is of the form $(\sum_{i=1}^N c_i s_i)^n$, where $c_1,...,c_N\in \K$, and, therefore, it is a value of the function
$$\K^n\ni (x_1,...,x_N)\to (\sum_1^N x_i s_i)^n\in \K\S,$$
which can be thought as a polynomial in $\S[x_1,...,x_N]$.
Hence, $I$ is generated by coefficients of such polynomials,
\be{equation}\label{sum}
\sum_{\sigma\in S_k} w_{\sigma(1)}^{\alpha_1}... w_{\sigma(k)}^{\alpha_k},
\en{equation}
for all sequences $\alpha_1\geq \alpha_2 \geq ...\geq \alpha_k\geq 1$, $1\leq k\leq n$, such that $\sum_{i=1}^k \alpha_i=n$ and for all $k$-tuples of distinct elements $w_1,...,w_k$ in $\S$.
Since the above generators of $I$ are homogeneous, the statement follows.
\epr

Since $I$ is generated by homogeneous elements of degree $\geq n$, $I\cap \K\S_d=0$ for $d<n$. Hence, $\D$ always contains all elements of degree less than $n.$ The elements of degree $d\geq n$ in $\D$ can be determined immediately from Lemma \ref{algor} by finding a basis of $\K\S_d/I\cap \K\S_d$.

In practice, the above approach works for groups $\Gamma$ with reasonably small generating sets. For free groups of large rank one can use the following method of finding generators:


\be{pro}
Suppose that $\B$ is a set of words defined in Theorem \ref{sln-main}(2) for the free semigroup of rank $\nu_n-1$, $\S=\langle s_1,...,s_{\nu_n-1}\rangle$ .
Consider a semigroup $\S'=\langle s_1',...,s_{N}'\rangle$, for $N\geq \nu_n-1,$
and the set $\B'\subset \S'$ composed of all words $s_{f(i_1)}'...s_{f(i_k)}'$ obtained by applying all strictly increasing functions $f:\{1,...,k\}\to \{1,...,N\}$
to all words $s_{i_1}...s_{i_k}\in \B.$  Then $\B'$ satisfies the condition of Theorem \ref{sln-main}(2) and, hence, it yields a generating set of $\K[\X_{SL(n,\K)}(\Gamma)]$, for the free group $\Gamma$ on $N$ generators.
\en{pro}

\bpr
Let $\D'\subset \S'$ be the set obtained by applying all strictly increasing functions $f:\{1,...,k\}\to \{1,...,N\}$ to the indices of words $s_{i_1}...s_{i_k}$ in $\D.$
Observe that $\D'\cap \S_{d,r}'$ spans $\K\S_{d,r}'/I\cap \K\S_{d,r}'.$ Indeed, suppose that $w\in \S_{d,r}'$ represents a non-zero element of $\K\S_{d,r}'/I\cap \K\S_{d,r}'$. Since $w$ is of length $k\leq \nu_n-1$,
$w$ is obtained by applying an increasing function $f:\{1,...,k\}\to \{1,...,N\}$ to the indices of the letters in a word $v\in \S_{d,r}$. However, $v$ (considered as an element of $\K\S_{d,r}/I\cap \K\S_{d,r}$) is a linear combination of words in $\D\cap \S_{d,r}.$  By applying $f$ to the indices of the letters in the components of this linear combination, we conclude that $w$ is a linear combination of words in $\D'$.

Now that we know that $\D_r'\cap \S_{d,r}'$ spans $\K\S_d'/I\cap \K\S_d',$ it is easy to see that
$\B'$ described above coincides with the one defined by Theorem \ref{sln-main}(2):
$$\B'=\{s_1',...,s_N'\}\cup \bigcup_{r=1}^N \D_r \cdot s_r'.$$
\epr


%
\section{Efficient generators for $Sp(n,\K)$- and $SO(n,\K)$-character varieties}
\label{s_spso}
%

In this section we formulate versions of Theorem \ref{sln-main} for symplectic and orthogonal groups.

We say that $*:\S\to \S$ is an involution on a semigroup $\S$ iff
$s^{**}=s$ and $(st)^*=t^*s^*$ for all $s,t\in \S.$
We say that $(\S,*)$ is a graded semigroup with an involution $*$ if $\S$ is graded, $*$ is an involution on $\S$, and $deg(s^*)=deg(s)$ for every $s\in \S.$

Let $\S_d=\{s\in \S: deg(s)\leq d\}.$
For $n$ even, let $I^s$ be the two sided ideal in $\K\S$ generated by the $n/2$-th powers of all symmetric elements, $(x+x^*)^{n/2},$ $x\in \K\S.$
Let $\D^s\subset \S$ be such that $\D^s\cap \S_d$ spans $\K\S_d/I^s \cap \K\S_d$ for all $d\in \N.$ (The upper index ``$s$" here stands for ``symmetric".)

Every group $\Gamma$ is a semigroup with an involution $\gamma^*=\gamma^{-1}.$

\begin{theorem}[Proof in sec \ref{s_sp-main}]\label{sp-main}\ \\
Let $\S$ be a graded semigroup with an involution. Assume that $\S$ is generated (as a semigroup with an involution) by elements $s_1,...,s_N$ of degree $1$.\\
\textup{(1)} For every epimorphism of semigroups with an involution $\eta: \S\to \Gamma,$
the $\K$-algebra $\K[\X_{Sp(n,\K)}(\Gamma)]$ is generated by $\sigma_{\gamma},$
for $\gamma\in \eta(\B^s),$ where
$$\B^s=\{s_1,...,s_N\}\cup \D^s \cdot \{s_1,...,s_N\},$$
for any $\D^s$ as above.

\textup{(2)} (Stronger version) Let $\S_{d,r}$ be the set of all elements of $\S$ which are represented by words in $\gamma_1,...,\gamma_r, \gamma_1^*,...,\gamma_r^*$ of degree $\leq d.$
If $\D^s\cap \S_{d,r}$ spans \mbox{$\K\S_{d,r}/I \cap \K\S_{d,r}$} for all $1\leq r\leq N$
and all $d\geq 1$, then $\K[\X_{Sp(n,\K)}(\Gamma)]$ is generated by
$$\B^s=\{s_1,...,s_N\}\cup \bigcup_{r=1}^N \D^s_r \cdot s_r,$$
where $\D_r^s\subset \D^s$ is the set of those elements of $\D^s$ which are words in $s_1,...,s_r,s_1^*,...,s_r^*$.
\end{theorem}

The algorithm of Sec. \ref{s_sln} can be easily modified to give an explicit method for finding $\B^s$.


%
%

Let $G=SO(n,\K)$ now.
By Corollary \ref{trace-epi}, any generating set of the $SL(n,\K)$-character variety yields a generating set of the $SO(n,\K)$-character variety for $n$ odd.
For $n$ even, the description of generators is more complicated, since (as noted in Sec. \ref{s_main}), the $G$-trace algebra of $\Gamma$ is a proper subring of $\K[\X_G(\Gamma)].$

Let $\B$ be as above. Let ${\cal M}\subset \B$ be such that $\{s-s^*: s\in {\cal M}\}$ spans the space
$\{s-s^*: s\in \K\B\}.$ In particular, one can take $\cal M$ to be the subset
of $\B$ composed of elements $s_{i_1}'....s_{i_d}'$ where $s_i'$ is either $s_i$ or
$s_i^*$ and the number of stared letters is not larger than the number of the non-stared ones.

\begin{theorem}[Proof in Sec. \ref{s_so2n-main}]\label{so2n-main}
Let $\S$ be a graded semigroup with an involution generated by elements $s_1,...,s_N$ of degree $1$.
Let $n$ be even. For every epimorphism of semigroups with an involution $\eta: \S\to \Gamma,$  $\K[\X_{SO(n,\K)}(\Gamma)]$ is generated as an $\T_{SO(n,\K)}(\Gamma)$-algebra by
$Q_n(w_1,...,w_{n/2}),$ for all possible $w_1,...,w_{n/2}\in \eta({\cal M}).$ ($Q_n$ was defined by (\ref{def-Q})).
\end{theorem}

Since $\K\S/I$ is spanned by monomials in $s_1,...,s_N,s_1^*,...,s_N^*$ of degree at most $\nu_n-1,$ Theorem \ref{basic_SO(2m)} follows.

%
\section{Proof of Proposition \ref{SL2-free}}
\label{s_proof_SL2}
%


(1) Let $\S$ be the free semigroup on $s_1,...,s_{N}$ and let
$\eta: \S\to \Gamma$ send $s_1,...,s_N$ to $\gamma_1,...,\gamma_N$.
By (\ref{e-Ku}), the nilpotency index of $\K\S/I$ is $3$ for $n=2.$ Hence, we can take $\D=\{s_1,...,s_{N}, s_is_j, 1\leq i,j\leq {N}\}.$
$$s_is_j+s_js_i=s_i^2+s_is_j+s_js_i+s_j^2=
(s_i+s_j)^2=0$$
implies $s_js_i=-s_is_j$ in $\K\S/I$.
Therefore, we can assume $i<j$ for $s_is_j$ in $\D$ and, hence,
by Theorem \ref{sln-main}(2),
$$\B=\{s_i,\ \text{for}\ i=1,...,N,\ s_is_j,\ \text{for}\ i\leq j,\
s_is_js_k,\ \text{for}\ i<j\leq k\}.$$
Since $\tau_{\alpha \beta^2}=\tau_\beta\tau_{\alpha \beta}+\tau_\alpha$ for all $\alpha,\beta\in \Gamma,$ we can always assume that $j<k$ above.

(2) Let $\S$ and $\eta:\S\to \Gamma$ be as above.
Since $\nu(\K\S/I)=2$, we take $\B=\{s_1,...,s_{N}\}$ and get
$$\B=\{s_i, i=1,..,N, s_is_j,\ \text{for}\ i\leq j\}.$$
Since $\tau_{\gamma^2}=\tau_\gamma^2-2$ we can assume that $i<j$ and the statement follows.
\qed

%
\section{Proof of Theorem \ref{normalization} and Proposition \ref{E_6}}
\label{s_norm}
%

Let $\Gamma$ be the free group on $N$ generators.\\
{\bf Proof of Theorem \ref{normalization}(1):} Let us assume that $N>1$ first.
The embedding $\rho(G)\subset SL(n,\K)$ induces a map
$$f: \rho(G)^{N}//{\cal N}(\rho(G))\to SL(n,\K)^N//SL(n,\K),$$
where ${\cal N}(\rho(G))$ is the normalizer of $\rho(G)$ in $SL(n,\K)$. (One can prove that ${\cal N}(\rho(G))$ is reductive and, hence, the quotient exists.)
By \cite{Vi}, $Im\, f$ is a closed subset of $SL(n,\K)^N//G.$ By the definition of the trace algebra, we have $\T_{\rho(G)}(\Gamma)= \K[Im\, f]$. By \cite{Vi}, $f$ is a normalization map onto its image.
Therefore, to complete (1) one needs to prove that
$$G^N//G \to \rho(G)^N//\rho(G) \to \rho(G)^N//{\cal N}(\rho(G))$$
is an isomorphism if $\rho$ is faithful and all twins of $\rho$ are identical.

The first map, $G^N//G \to \rho(G)^N//\rho(G)$, is an isomorphism since it is a bijection and $\rho(G)^N//\rho(G)$ is normal, cf. \cite[Sec 4.8]{Do}.
We claim that the map $\rho(G)^N//\rho(G) \to \rho(G)^N//{\cal N}(\rho(G))$ is an isomorphism as well:
The action of any element of ${\cal N}(G)$ by conjugation is an automorphism of $G$.
Therefore, if $\phi, \phi'$ had the same image in $\rho(G)^N//{\cal N}(\rho(G))$, then they would be related by an automorphism of $G$ and, hence, by our assumptions, $[\phi]= [\phi']$ in $\rho(G)^N//\rho(G)=G^N//G.$

For $N=1$, it is enough to show Vinberg's theorem for $N=1$, i.e. that
$$f: \rho(G)//{\cal N}(\rho(G))\to SL(n,\K)//SL(n,\K)$$
is a normalization map onto its image. By Example \ref{Z}, it is enough to prove that
$$T/W\to \rho(G)//{\cal N}(\rho(G))\to SL(n,\K)//SL(n,\K)$$
is a normalization map, where $T$ and $W$ are the maximal torus and the Weyl group of $G.$ That follows from the proof of Vinberg's result, since $T$ is generated by a single generic element.
(The only reason Vinberg needs $N>1$ is to make sure that his group $H$ is generated by $N$ generic elements.)

{\bf Proof of Theorem \ref{normalization}(2):} If a non-trivial central element $c$ belongs to $Ker\, \rho$, then
$[(c,...,c)]\ne [(e,...,e)]$ in $G^N//G$ have the same image under
$\pi: G^N//G\to \rho(G)^N//{\cal N}(\rho(G))$.
Since $\rho(G)^N//{\cal N}(\rho(G))$ is normal, cf. \cite[Sec 4.8]{Do}, and $\pi$ is not $1$-$1$,
the composition of $\pi$ with $f$ cannot be a normalization map.
\qed

{\bf Proof of Proposition \ref{E_6}:} As above, it is enough to show that the map
$\rho(G)^N//\rho(G) \to \rho(G)^N//{\cal N}(\rho(G))$ is an isomorphism.
The action of any element of ${\cal N}(\rho(G))$ by conjugation is an automorphism of $\rho(G)$.
It cannot be a non-inner automorphism of $G=E_6$ since any such automorphism sends $\rho$ to an inequivalent representation. Hence ${\cal N}(\rho(G))$ acts by inner automorphisms on $\rho(G)$ and the statement follows.\qed

%
\section{Proof of Theorem \ref{sln-main}}
\label{s_sln-main}
%

The proof is inspired by the work of Procesi, \cite{Pro-p,Pro-b}.
Recall that a $\K$-algebra $R$ is graded if $R=\bigoplus_{k=0}^\infty R_k$
as a vector space and $R_k\cdot R_l\subset R_{k+l}.$
An element of $R$ is homogeneous if it belongs to $R_k$ for some $k.$
All graded algebras in this paper are connected, i.e. $R_0=\K.$
An element $r\in R$ has degree $d$ if $d$ is the smallest index such that
$r\in \bigoplus_{k=0}^d R_k.$
We denote $\bigoplus_{k>0} R_k$ by $R^+.$

\be{lem}\label{graded}
If $R$ is a graded ring, then every element of degree $d$ in $R^+\cdot R^+$ belongs to the subring of $R$ generated by the homogeneous elements of $R$ of degree $<d.$
\en{lem}

\bpr
By splitting elements of $R^+$ into sums of homogeneous elements,
every element of $R^+\cdot R^+$
can be written as $r=\sum_i s_i\cdot t_i,$ where all $s_i,t_i$ are homogeneous of positive degree. After eliminating all summands such that
$deg\, s_i+deg\, t_i> deg\, r,$
the equality $r=\sum_i s_i\cdot t_i$ still holds and
$deg\, s_i,deg\, t_i<deg\, r$ for all $i.$
\epr

\be{rem}\label{Reynolds} Throughout the paper we will often use the following fact: If a reductive group $G$ acts on $\K$-algebras $A$ and $B$ such that an epimorphism $\phi: A\to B$ is $G$-equivariant, then $\phi$ restricts to a $G$-equivariant epimorphism
$A^G\to B^G.$ This follows from the complete reducibility of representations of reductive groups or, equivalently, from the properties of Reynolds operators.
\en{rem}

Let $C_{n,N}=\K[a_{i,j,k}:\, i=1,...,N,\, j,k=1,...,n].$
(The letter ``$C$" is used here to indicate that this will be our ring of coefficients.)
Let $FSG(s_1,...,s_N)$ be the free semigroup on $s_1,...,s_N$. It is a graded semigroup with $s_1,...,s_N$ having degree $1$. By the assumptions of Theorem \ref{sln-main}, there is an epimorphism of graded semigroups $FSG(s_1,...,s_N)\to \S$. By abuse of notation, the images of $s_1,...,s_N$ in $\S$ are denoted by the same symbols.
By an ``abstract nonsense" argument, there exists a unique universal quotient, $R(\S),$ of $C_{n,N}$ such that the homomorphism of semigroups
\be{equation}\label{Psi}
\Psi:FSG(s_1,...,s_N)\to M(n,C_{n,N})
\en{equation}
sending $s_1,...,s_N$ to
$$A_1=(a_{1,j,k})_{j,k=1,..,n},\ ...,\ A_N=(a_{N,j,k})_{j,k=1,..,n}\in M(n,C_{n,N})$$
composed with the natural projection $M(n,C_{n,N})\to M(n,R(\S))$
factors to a homomorphism
\be{equation}\label{P_NtoM_n}
\Psi: \S\to M(n,R(\S)),
\en{equation}
which, by abuse of notation, we denote by the same letter, $\Psi,$ as (\ref{Psi}).
(``Universal" means that every other such quotient factors through this one.) We call $R(\S)$ the universal representation $\K$-algebra in dimension $n$ and $\Psi$ the universal $n$-dimensional representation of $\S.$ This construction is analogous to that for groups
mentioned in Sec. \ref{s_sln}.

Here is a concrete construction of $R(\S):$ The semigroup $\S$ has a presentation
$$\S=\la s_1,...,s_N\ |\ r_{1,i}=r_{2,i},\ i\in {\cal I}\ra.$$
Then $R(\S)$ is the quotient of $C_{n,N}$ by the ideal generated by
$n^2$ entries of the matrix $\Psi(r_{1,i})-\Psi(r_{2,i})$
taken for every $i\in {\cal I}.$

Since $s_1,...,s_N$ are generators of $\S$ of degree $1$,
$deg(r_{1,i})=deg(r_{2,i})$ for each $i.$ Therefore, the grading on $C_{n,N}$
in which all generators $a_{i,j,k}$ have degree $1$ descends to a grading on $R(\S)$ (in which the image of every $a_{i,j,k}$ in $R(\S)$ has degree $1$). This grading will be important later.

Let $G=SL(n,\K).$
By the construction of $R(\S)$, there is a natural epimorphism $R(\S)\to R(\Gamma,G)$, where
$R(\Gamma,G)$ is the universal representation algebra of Sec. \ref{s_sln}.

Since $C_{n,N}$ is the coordinate ring of $M(n,\K)^N$ (the Cartesian product of $N$ copies of $M(n,\K)$), the $SL(n,\K)$ action on $M(n,\K)^N$ by conjugation induces an $SL(n,\K)$ action on $C_{n,N}.$
This action descends to an action on $R(\S)$ and on $R(\Gamma,G)$.
Hence, we have an epimorphism
$$R(\S)^G\to R(\Gamma,G)^G=\K[\X_G(\Gamma)]$$
by Remark \ref{Reynolds}.

By abuse of notation, denote the images of $A_1,...,A_N$ under the projection $M(n,C_{n,N})\to M(n,R(\S))$ by the same symbols.
Let $T(\S)$ be the subalgebra of $R(\S)$ generated by the traces of monomials
in $A_1,...,A_N\in M(n,R(\S))$.
Clearly, $T(\S)\subset R(\S)^G.$

\be{lem}
$T(\S)= R(\S)^G.$
\en{lem}

\bpr
By \cite{Pro-p,Pro-b}, $C_{n,N}^G$ is generated by the traces of monomials
in $A_1,...,A_N\in M(n,C_{n,N}).$ Now the statement follows from the fact that
$R(\S)$ is a $G$-equivariant quotient of $C_{n,N}$ and from Remark \ref{Reynolds}.
\epr

Let $\D,\B \subset \S$ be defined as in the statement of Theorem \ref{sln-main}.
By abuse of notation, denote the set of monomials in $A_1,...,A_N$ in
$M(n,R(\S))$ corresponding to the elements of $\D$ and of $\B$ in $\S$ via the map (\ref{P_NtoM_n}) by the same symbols, $\D,\B.$

We are going to complete the proof of Theorem \ref{sln-main} (stronger version), by showing that $T(\S)$ is generated by traces of monomials in $\B.$
We will use induction on the degree.

\begin{pro}\label{sln-red}
For every monomial $M$ of degree $d>1$ in variables $A_1,...,A_N,$  $tr(M)$ belongs to the $\K$-subalgebra of $R(\S)^G$ generated by the traces of monomials in $\Psi(\B)$ and by the traces of monomials in $A_1,...,A_N$ of degree $<d.$
\end{pro}

\noindent{\it Proof.}
Identify $R(\S)$ with the scalar matrices in $M(n,R(\S)).$
Let $S(\S)$ be the $T(\S)$-subalgebra of $M(n,R(\S))$ generated by the matrices $A_1,...,A_N.$

Note that $T(\S)$ is a graded subalgebra of $R(\S).$
Additionally $M(n,R(\S))$ is a graded algebra, with a matrix being a homogeneous element of degree $k$ iff all its entries are homogeneous of degree $k$ in $R(\S).$ In particular,
$A_{i_1}...A_{i_k}$ is homogeneous of degree $k.$
$S(\S)$ is a graded subalgebra of $M(n,R(\S)).$

Define $T(\S)^+\subset T(\S)$ and $S(\S)^+\subset S(\S)$
to be the subalgebras without identity spanned by the homogeneous elements of positive degree.
Hence,
$$T(\S)=T(\S)^+\oplus \K\ \text{and}\ S(\S)=S(\S)^+\oplus \K.$$

The homomorphism (\ref{P_NtoM_n}) extends to
a homomorphism
$$\Psi: \K\S \to S(\S)^+$$
sending $s_i$ to $A_i.$

\be{lem}\label{ITS}
$$\Psi(I)\subset T(\S)^+S(\S).$$
\en{lem}

\bpr
Every element of $\Psi(I)$ is a sum of elements $XY^nZ$, where $X,Z\in S(\S)$ and $Y\in S(\S)^+.$
The matrix $Y\in S(\S)^+$ satisfies its characteristic polynomial
$$Y^n+\sum_{i=0}^{n-1} c_iY^i=0,$$
with $c_0,...,c_{n-1}\in T(\S)$, since they are conjugation invariant.
Furthermore, they belong to $T(\S)^+$ since each $c_i$ is homogeneous of degree $n-i$ in the entries of $Y\in S(\S)^+.$
Hence,
$$XY^nZ=-\sum_{i=0}^{n-1} c_iXY^iZ\in T(\S)^+S(\S).$$
\epr

Let $M=A_{i_1}....A_{i_d},$ for some $i_1,..,i_d\in \{1,...,N\},$ $d>1.$
Let $r= Max\{i_1,...,i_d\}.$ Since $Tr(A_{i_1}...A_{i_d})$ is invariant under a cyclic permutation of $A_{i_1},....,A_{i_d},$ we can assume that $i_d=r.$
By the definition of $\D$,
\begin{equation}\label{reduction-sl}
s_{i_1}...s_{i_{d-1}}=\sum_{s\in \D_r,\, deg\, s\leq d-1}\, c_s\cdot s +C,
\end{equation}
where $c_s\in \K$ and $C\in I.$

By multiplying both sides of (\ref{reduction-sl}) by $s_{i_d}$,
applying $\Psi$ and then taking trace, we get
$$Tr(M)=Tr(A_{i_1}...A_{i_d})= \sum_{s\in \D_r,\, deg\, s\leq d-1}\, c_s\cdot Tr(\Psi(s) A_r) + Tr(\Psi(C)A_r).$$

Note that $\Psi(s) A_r\in \B$ for $s\in \D_r$
and that $Tr(\Psi(C) A_r)$ is an element of degree $\leq d$ in $T(\S).$ Furthermore, by Lemma \ref{ITS},\vspace*{.1in}\\
$Tr(\Psi(C)A_r)\in Tr(T(\S)^+S(\S)A_r)=Tr(T(\S)^+S(\S)^+)=$\hspace*{1in}\\
\ \hspace*{2in} $T(\S)^+\cdot Tr(S(\S)^+)=T(\S)^+\cdot T(\S)^+.$\vspace*{.1in}

Since the traces of monomials in $A_1,...,A_N$ are homogeneous generators of $T(\S),$ Lemma \ref{graded} implies that $Tr(\Psi(C)A_r)$ is polynomial in traces of monomials of degree $<d.$ This completes the proof of Proposition \ref{sln-red}.
\qed

%
\section{Proof of Theorem \ref{spso}}
\label{spso-proof}
%

We are going to prove a stronger statement: For every $\Gamma$ and every $G$ symplectic, orthogonal, or special orthogonal, $\K[\X_G(\Gamma)]$ is generated by $\sigma_\gamma,$ for $\gamma\in \Gamma.$

We have a natural epimorphism $C_{n,N}\to R(\Gamma,G)$, as in the previous section, inducing a homomorphism
$M(n, C_{n,N})\to M(n, R(\Gamma,G))$ sending $A_i=(a_{i,j,k})_{j,k=1,...,n}$ for $i=1,...,N,$ to the images of
$\gamma_1,...,\gamma_N$ under the universal representation $\Psi: \Gamma\to M(n, R(\Gamma,G)).$

The $G$ action on $M(n,\K)^N$ by conjugation defines a $G$ action on $C_{n,N}$,
which descends to an action on $Hom(\Gamma,G).$
By Remark \ref{Reynolds},
\be{equation}\label{spso-map}
C_{n,N}^G\to R(\Gamma,G)^G=\K[\X_G(\Gamma)].
\en{equation}
is an epimorphism.

By \cite[Thm 10.1]{Pro-p}, \cite[Sec 11.8.2]{Pro-b}, $C_{n,N}^G$ is generated by
the traces of monomials in the matrices $A_1,...,A_N,A_1^*,...,A_N^*,$ where
$A^*=JA^TJ^{-1},$ for $G$ symplectic and $A^*=A^T,$ for $G$ orthogonal and special odd orthogonal.
Since (\ref{spso-map}) maps each such trace to $\sigma_\gamma$ for $\gamma$ being the corresponding word in $\gamma_1,...,\gamma_N,\gamma_1^*,...,\gamma_N^*$, the statement follows.
\qed

%
\section{Proof of Theorem \ref{sp-main}}
\label{s_sp-main}
%

The proof is an adaptation of the proof of Theorem \ref{sln-main}.

As before, let $FSG(s_1,...,s_N,s_1^*,...,s_N^*)$ be the free semigroup on
$s_1,...,s_N,$ $s_1^*,...,s_N^*$ (and, hence, the free semigroup with an involution on $s_1,...,s_N$). Consider $M(n,C_{n,N})$ as a semigroup with the symplectic involution, $A^*=JA^TJ^{-1}.$

Let
$$\Psi: FSG(s_1,...,s_N,s_1^*,...,s_N^*)\to M(n,C_{n,N})$$
be the homomorphism of semigroups with involutions sending $s_i$ to $A_i=(a_{i,j,k})_{j,k=1,...,n}$
Let $R(\S)$ be the universal quotient of $C_{n,N}$ such that $\Psi$ composed with
$M(n,C_{n,N})\to M(n,R(\S))$ factors through
$$\Psi: \S\to M(n,R(\S)).$$
Since $\S$ has a homogeneous set of generators of degree $1$, the grading on $C_{n,N}$ descends to a grading on $R(\S)$.

As before, we have a natural epimorphism
$$R(\S)\to R(\Gamma,Sp(n,\K)).$$
The $Sp(n,\K)$ action on $M(n,\K)^N$ by conjugation defines an $Sp(n,\K)$ action on $C_{n,N}$,
which descends to an action on $R(\S)$ and on $R(\Gamma,Sp(n,\K)).$
By Remark \ref{Reynolds},
$$R(\S)^{Sp(n,\K)}\to R(\Gamma,Sp(n,\K))^{Sp(n,\K)}=\K[\X_{Sp(n,\K)}(\Gamma)].$$
is an epimorphism.

By abuse of notation, denote the images of $A_1,...,A_N, A_1^*,...,A_N^*\in M(n,C_{n,N})$ in $M(n,R(\S))$ (under the natural projection) by the same symbols.
Let $T(\S)$ be the subalgebra of $R(\S)$ generated by the traces of monomials
in these matrices. Then $T(\S)= R(\S)^{Sp(n,\K)},$ by \cite[Thm 10.1]{Pro-p}, \cite[Sec 11.8.2]{Pro-b}.

Therefore, in order to establish Theorem \ref{sp-main} (stronger version) it is enough to prove that $T(\S)$ is generated by traces of monomials in $\B^s$.
That will follow by induction on the degree of monomials, from the following proposition.
(The degree of a monomial is the number of components $A_i$ and $A_i^*$
in it.)

\begin{pro}\label{sp-red}
For every monomial $M$ of degree $d>1$ in variables $A_1,...,A_N,$ $A_1^*,...,A_N^*$, $tr(M)$ belongs to the $\K$-subalgebra of $T(\S)$ generated by the traces of monomials in $\Psi(\D^s)$ and by the traces of monomials of degree $<d.$
\end{pro}

The proof of this proposition is almost identical to that of Proposition \ref{sln-red}:
Let $M=A_{i_1}'....A_{i_d}',$ for some $i_1,..,i_d\in \{1,...,N\},$ $d>1,$
where $A_i'$ is either $A_i$ or $A_i^*.$
Let $r= Max\{i_1,...,i_d\}.$ Since $tr(M)=tr(M^*),$
one can replace $M$ with $M^*$ without loss of generality. Therefore, we can assume that there is at least one component $A_r$ (without the star) in that monomial.
Since $Tr(A_{i_1}...A_{i_d})$ is invariant under a cyclic permutation of its components, we can assume that $i_d=r.$
By the definition of $\D^s$,
\begin{equation}\label{reduction-sp}
s_{i_1}...s_{i_{d-1}}=\sum_{s\in \D^s_r,\, deg\, s\leq d-1}\, c_s\cdot s +C,
\end{equation}
where $c_s\in \K$ and $C\in I^s.$

By multiplying both sides of (\ref{reduction-sp}) by $s_{i_d}$,
applying $\Psi$ and then taking trace, we get
$$Tr(M)=Tr(A_{i_1}'...A_{i_d}')= \sum_{s\in \D,\, deg\, s\leq d-1}\, c_s\cdot Tr(\Psi(s) A_r) + Tr(\Psi(C)A_r).$$
The completion of the proof is identical to that of Proposition \ref{sln-red}, except that
one replaces Lemma \ref{ITS} with the following one:

\be{lem}\label{ISTS}
$$\Psi(I^s)\subset T(\S)^+S(\S).$$
\en{lem}

\bpr
Every element of $\Psi(I^s)$ is a sum of elements $X(Y+Y^*)^{n/2}Z$, where $X,Z\in S(\S)$ and $Y\in S(\S)^{+}.$
For $n$ even and $M\in M(n,\K)$ invariant under the symplectic involution,
$$Pf_M(\lambda)=\pf((\lambda I-M)J)$$
is called the characteristic Pfaffian of $M.$
(Here $I$ is the identity matrix and $J$ the skew-symmetric matrix used to define $Sp(n,\K)$.)

Let
$$Pf_M(\lambda)=\lambda^n+\sum_{i=0}^{n-1} c_i\lambda^i$$
be the characteristic Pfaffian of $M=X+X^*$.
Since every matrix invariant under the symplectic involution satisfies its characteristic Pfaffian equation, $$Pf_M(M)=0,$$ cf. \cite[Sec 11.8.7]{Pro-b}.
The coefficients $c_0,...,c_{n-1}$ belong to $R(\S)^G=T(\S)$, since they are conjugation invariant.
Furthermore, they belong to $T(\S)^+$ since each $c_i$ is homogeneous of degree $n-i$ in the entries of $M\in S(\S)^+.$

Hence,
$$X(Y+Y^*)^nZ=-\sum_{i=0}^{n-1} c_iX(Y+Y^*)^iZ\in T(\S)^+S(\S).$$
\epr


%
\section{Proof of Theorem \ref{so2n-main}}
\label{s_so2n-main}
%

Let $G=SO(n)$ for $n$ even. Let $R(\S)$ be defined as in the proof of Theorem \ref{sp-main}.
Let $T(\S)$ be the subalgebra of $R(\S)$ generated by the traces of monomials in matrices $A_1,...,A_N,A_1^*,...,A_N^*\in M(n,R(\S)).$

Since there is the natural epimorphism
$$R(\S)^{G}\to R(\Gamma,G)^{G}=\K[\X_{G}(\Gamma)],$$
mapping $T(\S)$ onto $\T_G(\Gamma)$, the statement of Theorem \ref{so2n-main} is implied by the following:

\be{pro}
\textup{(1)} $R(\S)^G$ is generated by $T(\S)$ and by the values of\\ $Q_n(M_1,...,M_{n/2})$
for all $M_1,...,M_{n/2}$ in $S(\S).$\\
\textup{(2)} It is enough to consider values of $Q_n$ for $M_1,...,M_{n/2}\in \D$ only.\\
\textup{(3)} Furthermore, it is enough to consider $M_1,...,M_{n/2}\in {\cal M}$ only.
\en{pro}

\bpr (1) By Remark \ref{Reynolds}, the natural projection $M(n,C_{n,N})^G\to R(\S)^G$ is onto. Therefore, it is enough to prove that statement for the semigroup $M(n,C_{n,N})$ with involution $M^*=M^T$. That was done in \cite[Sec 11.8.2]{Pro-b}.\\
We will prove (2) by contradiction: Denote by $T'(\S)$ the subalgebra of
$R(\S)^G$ generated by $T(\S)$ and by the values of $Q_n(M_1,...,M_{n/2})$
for $M_1,...,M_{n/2}\in \D$ only. Assume that $T'(\S)\ne T(\S).$
Let $(d_1,...,d_{n/2})\in \N^{n/2}$ be the smallest element in the lexicographic order such that there exist monomials $M_1,...,M_{n/2}$ of degrees $d_1,...,d_{n/2}$ such that $Q_n(M_1,...,M_{n/2})\not\in T'(\S).$
Then at least one $M_i$ is not in $\Psi(\D).$ Let us assume that it is the first one for simplicity. Abbreviate $d_1$ to $d.$
Then $M_1=A_{i_1}'...A_{i_d}',$ where each $A_i'$ is $A_i$ or $A_i^*.$
By the definition of $\D,$
$$s_{i_1}'...s_{i_d}'=\sum_{s\in \D,\ deg\, s\leq d} c_s\cdot s+ C,$$
where $C\in I.$

By applying $\Psi$ we get
\be{equation}\label{MSC}
M_1=\sum_{s\in \D,\ deg\, s\leq d} c_s\cdot \Psi(s)+ \Psi(C).
\en{equation}
By Lemma \ref{ITS},
\be{equation}\label{CXY}
\Psi(C)=\sum_j T_j\cdot X_j,
\en{equation}
where $T_j\in T^+(\S),$ $X_j\in S(\S).$ Furthermore, by splitting each $T_j,X_j$
into a sum of homogeneous summands if necessary, we can assume that all
$T_j,X_j$ are homogeneous with respect to the grading on $T(\S)$ and on $S(\S).$
Finally, we can remove all summands $T_jX_j$ of degree greater than $d$ from the right side of (\ref{CXY}) without loss of validity of that equation.
Therefore, (\ref{MSC}) combined with (\ref{CXY}) expresses $M_1$ as a sum of terms
of the form $TX$ where $T\in T(\S)$, $X\in S(\S)$ and and either $X\in \Psi(\D)$ or $deg\, X<d.$
Since elements of $T(\S)$ are scalar matrices in $S(\S)$ and
$Q_n$ is multi-linear,
$Q_n(M_1,...,M_{n/2})$ can be expressed as a sum of terms $T\cdot Q_n(X,M_2,...,M_{n/2})$ with $deg\, X<deg\, M_1.$ Therefore,
$Q_n(X,M_2,...,M_{n/2})\not\in T'(\S),$ for at least one such $X$ of degree $<d$ --
contradicting the initial assumption.

(3)
Since the substitution of $M_i$ by $M_i^T$ changes sign in $Q_n(M_1,...,M_{n/2}),$
it is enough to consider $M_1,...,M_{n/2}$ in the skew-symmetric part of $\K\D$ only. Elements $\{M-M^T: M\in \Psi({\cal M})\}$ span that space.
Finally, since the substitution of $M_i-M_i^T$ by $M_i$ multiplies the value of $Q_n$ by $2,$ it is enough to consider $M_1,...,M_{n/2}\in {\cal M}.$
\epr

\centerline{Dept. of Mathematics, 244 Math. Bldg.}
\centerline{University at Buffalo, SUNY}
\centerline{Buffalo, NY 14260, USA}
\centerline{asikora@buffalo.edu}

\end{document}